\documentstyle[12pt,leqno]{article}
\begin{document}
\title{The BV-algebra of a Jacobi manifold}
\author{Izu Vaisman}
\date{}
\maketitle
{\def\thefootnote{*}\footnotetext[1]%
{{\it 1991 Mathematics Subject Classification}
17 B 66, 58 F 05. \newline\indent{\it Key words and phrases}:
Gerstenhaber algebras,
BV-algebras, Jacobi manifolds, Lie bialgebroids.}}
\begin{center} \begin{minipage}{12cm}
A{\footnotesize BSTRACT.
We show that the Gerstenhaber algebra of the $1$-jet Lie algebroid of a
Jacobi manifold has a canonical
exact generator, and discuss duality between its homology and
the Lie algebroid cohomology. We also give a new example of a Lie
bialgebroid on Poisson manifolds.}
\end{minipage} \end{center}
\section{Introduction}
A {\em Gerstenhaber algebra} is a triple $({\cal A}=\oplus_{k\in{\bf Z}}
{\cal A}^{k},\wedge,[\;,\;])$ where $\wedge$ is an associative, graded
commutative algebra structure (e.g., over {\bf R}), $[\;,\;]$ is a graded
Lie algebra structure for the {\em shifted grades} $[k]:=k+1$ (the sign
$:=$ denotes a definition), and
$$[a,b\wedge c]=[a,b]\wedge c+(-1)^{kj}b\wedge[a,c], \leqno{(1.1)}$$
$\forall a\in{\cal A}^{k+1}$, $b\in {\cal A}^{j}$, $c\in{\cal A}$.
If this structure is supplemented by 
an endomorphism $\delta:{\cal A}
\rightarrow{\cal A}$, of grade $-1$,
such that $\delta^{2}=0$ and
$$[a,b]=(-1)^{k}(\delta(a\wedge b)-\delta a\wedge b-(-1)^{k}a\wedge\delta
b)\hspace{3mm}(a\in{\cal A}^{k},\,b\in{\cal A}),\leqno{(1.2)}$$
one gets an {\em exact Gerstenhaber algebra} or {\em Batalin-Vilkovisky
algebra} ({\em BV-algebra}) with the {\em exact generator} $\delta$.
If we also have a differential $d:{\cal A}^{k}\rightarrow
{\cal A}^{k+1}$ ($d^{2}=0)$ such that
$$d(a\wedge b)=(da)\wedge b+(-1)^{k}a\wedge(db) \hspace{3mm}
(a\in{\cal A}^{k},\,b\in{\cal A}),\leqno{(1.3)}$$
we will say that we have a {\em differential BV-algebra}.
(Some authors include $d$ in the BV-structure \cite{YK}, \cite{Xu}.)
Finally,if
$$d[a,b]=[da,b]+(-1)^{k}[a,db]\hspace{3mm}(a\in{\cal A}^{k},b\in{\cal
A})\leqno{(1.4)}$$
the differential BV-algebra is said to be strong \cite{Xu}.

On the other hand, a Jacobi manifold
(e.g., \cite{GL}) is a smooth manifold $M^{m}$
(everything is of class
$C^{\infty}$ in this paper) with a Lie algebra structure of
local type on the space of functions $C^{\infty}(M)$ or, equivalently
\cite{GL}, with
a bivector field $\Lambda$ and a vector field $E$ such that
$$[\Lambda,\Lambda]=2E\wedge\Lambda,\;\;[\Lambda,E]=0.\leqno{(1.5)}$$
In (1.5) one has usual
Schouten-Nijenhuis brackets. If $E=0$, $(M,\Lambda)$
is a {\em Poisson manifold}.

One of the most interesting examples of a BV-algebra is that of the
Gerstenhaber algebra of the cotangent Lie algebroid of a Poisson manifold,
described by Koszul \cite{Kz} and
Y. Kosmann-Schwarzbach \cite{YK}. More generally, Xu \cite{Xu}
extends a result of Koszul \cite{Kz} and proves that the
exact generators of the Gerstenhaber
algebra of a Lie algebroid $A\rightarrow M$
are provided by flat connections on $\wedge^{r}A$
($r=rank\,A$), and Huebschmann \cite{H}
proves a corresponding result for {\em Lie-Rinehart algebras}.

The main aim of this note is to show that a Jacobi manifold also has a
canonically associated, differential,
BV-algebra (which, however, is not strong)
namely, the Gerstenhaber algebra of the $1$-jet Lie algebroid defined by
Kerbrat-Benhammadi \cite{KB}. Then, we apply results of Xu \cite{Xu},
and Evens-Lu-Weinstein \cite{ELW} 
to discuss duality between the homology of this
BV-algebra and the cohomology of the Lie algebroid. (The homology was also
independently introduced and studied by de Le\'on, Marrero and Padron
\cite{LMP0}.) 

In the final section, we come back to a Poisson
manifold $M$ with the Poisson bivector $Q$,
and show that the infinitesimal automorphisms $E$ of $Q$
yield natural Poisson bivectors of the Lie algebroid
$TM\oplus{\bf R}$. These bivectors lead to triangular Lie bialgebroids and
BV-algebras in the usual way \cite{YK}, \cite{Xu}

Notice that BV-algebras play an important role in some recent researches
of theoretical physics (e.g., \cite{G}).

{\it Acknowledgements}.
The final version of this paper was written during the author's visit at the
Centre de Math\'ematiques, \'Ecole Polytechnique, Palaiseau, France, 
and he wishes to thank his host institution for invitation
and support.
The author is grateful to Y. Kosmann-Schwarzbach
for her invitation to \'Ecole Polytechnique, and several useful
discussions on the modular class of a triangular Lie algebroid, as well as 
to J. C. Marrero for the comparison of the
BV-homology and that of \cite{LMP0}, and to J. Monterde for his careful
reading of the final text and his remarks.
\section{The Jacobi BV-algebra}
For any Lie algebroid
$A\rightarrow M$ with anchor $\alpha:A\rightarrow TM$
one has a Gerstenhaber algebra ${\cal A}(A)$ defined by
$${\cal A}(A):=(\oplus_{k\in{\bf
N}}\Gamma\wedge^{k}A,\wedge,[\;,\;]_{SN}),
\leqno{(2.1)}$$
where $\Gamma$ denotes spaces of global cross-sections, and
SN denotes the Schouten-Nijenhuis bracket
(e.g., \cite{YK}, \cite{Xu};
on the other hand, we refer the reader to \cite{{Mk},{YK},{ELW}},
for instance, for the basics of Lie algebroids and Lie algebroid
calculus).
The BV-algebra which we want to discuss is associated
with the $1$-jet Lie algebroid of a Jacobi manifold $(M,\Lambda,E)$ defined
in \cite{KB}, which we present as follows.

We identify $M$ with $M\times\{0\}\subseteq M\times{\bf R}$, where
$M\times{\bf R}$
is endowed with the Poisson bivector \cite{GL}
$$P:=e^{-t}(\Lambda+\frac{\partial}{\partial t}\wedge E)
\hspace{3mm}(t\in{\bf R}).\leqno{(2.2)}$$
Let $J^{1}M=T^{*}M\oplus{\bf R}$ be the vector bundle of $1$-jets of
real functions on $M$, and notice that $\Gamma J^{1}M$ is isomorphic as a
$C^{\infty}(M)$-module with
$$\Gamma_{0}(M):=\{e^{t}(\alpha+fdt)\,/\alpha\in\wedge^{1}M,\,f\in
C^{\infty}(M)\}\subseteq\wedge^{1}(M\times{\bf R}). \leqno{(2.3)}$$
A straightforward computation shows that $\Gamma_{0}(M)$ is closed 
under the
bracket of the cotangent Lie algebroid  of
$(M\times{\bf R},P)$ (e.g., \cite{V3}) namely,
$$\{e^{t}(\alpha+fdt),e^{t}(\beta+gdt)\}_{P}=
e^{t}[L_{\sharp_{\Lambda}\alpha}\beta-L_{\sharp_{\Lambda}\beta}\alpha
-d(\Lambda(\alpha,\beta)) \leqno{(2.4)}$$
$$+fL_{E}\beta-gL_{E}\alpha-\alpha(E)\beta+\beta(E)\alpha
+(\{f,g\}-\Lambda(df-\alpha,dg-\beta))dt],$$
where $<\sharp_{\Lambda}\alpha,\beta>:=\Lambda(\alpha,\beta)$
$(\alpha,\beta\in\wedge^{1}M)$, and
$$\{f,g\}=\Lambda(df,dg)+f(Eg)-g(Ef)\hspace{3mm}(f,g\in C^{\infty}(M))$$
is the bracket which defines the Jacobi structure \cite{GL}.

Therefore, (2.4) produces a Lie bracket on $\Gamma J^{1}M$. Moreover, if
$\sharp_{P}$ is defined similar to $\sharp_{\Lambda}$, we get
$$\sharp_{P}(e^{t}(\alpha+fdt))=\sharp_{\Lambda}\alpha+fE-\alpha(E)
\frac{\partial}{\partial t},\leqno{(2.5)}$$ and
$$\rho:=(pr_{TM}\circ\sharp_{P})_{t=0}:J^{1}M\rightarrow TM \leqno{(2.6)}$$
has the properties of an anchor, since so does $\sharp_{P}$.

Formulas (2.4), (2.6) precisely yield the Lie algebroid structure defined
in \cite{KB}. In what follows we refer to it as the {\em $1$-jet
Lie algebroid}. The mapping $f\mapsto e^{t}(df+fdt)$ is a Lie
algebra homomorphism from the Jacobi algebra of $M$ to $\Gamma_{0}M$.
\proclaim 2.1 Proposition. The Gerstenhaber algebra ${\cal A}(J^{1}M)$
is isomorphic to the subalgebra ${\cal A}_{0}(M):=
\oplus_{k\in {\bf N}}\wedge^{k}\Gamma_{0}(M)$ of the Gerstenhaber
algebra
${\cal A}(T^{*}(M\times{\bf R}))$.\par
\noindent{\bf Proof.} The elements of ${\cal A}^{k}_{0}(M)$
are of the form
$$\lambda=e^{kt}(\lambda_{1}+\lambda_{2}\wedge dt)\hspace{3mm}
(\lambda_{1}\in\wedge^{k}M,\,\lambda_{2}\in\wedge^{k-1}M),
\leqno{(2.7)}$$
and we see that ${\cal A}_{0}(M)$ is closed by the wedge product and by the
bracket $\{\;,\;\}_{P}$ of general differential forms on the Poisson manifold
$(M\times{\bf R},P)$ (e.g., \cite{V3}). Accordingly,
$({\cal A}(J^{1}M),\wedge,\{\:,\:\})$
and $({\cal A}_{0}(M),\wedge,\{\:,\:\}_{P})$ are isomorphic
Gerstenhaber algebras since they are isomorphic at the grade $1$ level,
and the brackets of terms of higher degree are spanned by those of
degree $1$.
Q.e.d.
\proclaim 2.2 Remark. Since ${\cal A}_{0}(M)$ is a Gerstenhaber algebra,
the pair $({\cal A}_{0}^{0}=C^{\infty}(M),\,{\cal
A}_{0}^{1}=\Gamma_{0}(M))$ is a Lie-Rinehart algebra \cite{H}.
\par
Now, we can prove
\proclaim 2.3 Proposition. The Gerstenhaber algebra ${\cal A}_{0}(M)$ has
a canonical exact generator. \par
\noindent{\bf Proof.} It is known that ${\cal A}(T^{*}(M\times {\bf R}))$
has the exact generator of Koszul and Brylinski (e.g., \cite{V3})
$$\delta_{P}=i(P)d-di(P),\leqno{(2.8)}$$
where $P$ is the bivector (2.2). Hence, all we have to do is check that
$\delta_{P}\lambda\in{\cal A}_{0}^{k-1}(M)$ if
$\lambda$ is given by (2.7).

First, we notice that
$$i(P)(dt\wedge\mu)=e^{-t}(i(E)\mu+dt\wedge(i(\Lambda)\mu))
\hspace{3mm}(\mu\in\Lambda^{*}M).\leqno{(2.9)}$$
Then, if we also introduce the operator
$\delta_{\Lambda}:=i(\Lambda)d-di(\Lambda)$ \cite{CLM}, and compute for
$\lambda$ of (2.7), we get
$$\delta_{P}\lambda=e^{(k-1)t}[\delta_{\Lambda}\lambda_{1}+(-1)^{k}L_{E}
\lambda_{2}+ki(E)\lambda_{1} \leqno{(2.10)}$$
$$+(\delta_{\Lambda}\lambda_{2}+(-1)^{k}i(\Lambda)\lambda_{1}+
(k-1)i(E)\lambda_{2})\wedge dt].$$
Q.e.d.

Of course, $\delta_{P}$ of (2.10) translates to an exact generator $\delta$
of the Gerstenhaber algebra ${\cal A}(J^{1}M)$, and the latter becomes a
BV-algebra. This is the BV-algebra announced in Section 1, and we call
it the
{\em Jacobi BV-algebra of the Jacobi manifold}
$(M,\Lambda,E)$.
We can look at it under the two
isomorphic forms indicated by Proposition 2.1.

It is easy to see that the Jacobi BV-algebra above has the
differential $$\bar d\lambda:=e^{(k+1)t}d(e^{-kt}\lambda), \leqno{(2.11)}$$
where $\lambda$ is given by (2.7). But, $\bar d$ is not a derivation of the
Lie bracket $\{\;,\;\}$ of ${\cal A}(J^{1}M)$, and computations lead to
$$(\delta_{P}\bar d+\bar d\delta_{P})\lambda=e^{kt}[(k+1)i(E)d\lambda_{1}
\leqno{(2.12)}$$
$$+(L_{E}\lambda_{2}+(k+1)i(E)
d\lambda_{2}-(-1)^{k}\delta_{\Lambda}\lambda_{1})
\wedge dt],$$ where $\lambda$ is given by (2.7) again.
\proclaim 2.4 Remark. If we refer to the Poisson case $E=0$, we see that
both $T^{*}M$ and $J^{1}M$ have natural structures of Lie algebroids.
The Lie bracket and anchor map of $J^{1}M$ are given by
$$\{e^{t}(\alpha+fdt),e^{t}(\beta+gdt)\}=e^{t}[\{\alpha,\beta\}_{\Lambda}
+((\sharp_{\Lambda}\alpha)g-(\sharp_{\Lambda}\beta)f-\Lambda(\alpha,\beta)
)dt],\leqno{(2.13)}$$
$$\rho(e^{t}(\alpha+fdt))=\sharp_{\Lambda}\alpha, \leqno{(2.14)}$$
and the mapping $\alpha\mapsto e^{t}(\alpha+0dt)$ preserves the Lie
bracket, hence, $T^{*}M$ is a Lie subalgebra of $J^{1}M$,
and the latter is an extension of the former by the trivial line bundle
$M\times{\bf R}$. $J^{1}M$
was not yet used in Poisson geometry. \par
\section{The homology of the Jacobi BV-algebra}
We call the homology of the Jacobi BV-algebra of a Jacobi manifold
$(M,\Lambda,E)$, with boundary operator $\delta$, {\em Jacobi homology}
$H^{J}_{k}(M,\Lambda,E)$. (Another Jacobi homology was studied in \cite{CLM}.)
Here, we look at this homology from the point of view of
\cite{Xu} and \cite{ELW}, and discuss duality
between the Jacobi homology and the Lie algebroid cohomology of $J^{1}M$,
called {\em Jacobi cohomology}.

Jacobi cohomology coincides with the one
studied by de Le\'on, Marrero and Padr\'on in
\cite{LMP}. If $C\in\Gamma\wedge^{k}(J^{1}M)^{*}$ is seen as a $k$-multilinear
skew symmetric form on arguments (2.7) of degree $1$, at $t=0$, it
may be written as
$$C=\tilde C_{/t=0}:=e^{-kt}[(C_{1}+\frac{\partial}{\partial t}\wedge
C_{2})]_{t=0}
\hspace{3mm}(C_{1}\in{\cal V}^{k}M,\,C_{2}\in{\cal V}^{k-1}M),
\leqno{(3.1)}$$
where ${\cal V}^{k}M$ denotes the space of $k$-vector fields on $M$.
Furthermore, the coboundary, say $\sigma$, is given by the usual formula
$$(\sigma C)(s_{0},...,s_{k})=\sum_{i=0}^{k}(-1)^{i}(\rho s_{i})
C(s_{0},...,\hat s_{i},...,s_{k}) \leqno{(3.2)}$$
$$+
\sum_{i<j=0}^{k}(-1)^{i+j}C(\{s_{i},s_{j}\},s_{0},...,
\hat s_{i},...,\hat s_{j},...,s_{k}),$$
where $\rho$ is given by (2.6), and $s_{i}\in\Gamma J^{1}M$. Again, if we see
the arguments as forms (2.7) with $k=1$, 
(3.2) becomes
$$(\sigma C)=[\sigma_{P}\tilde C)]_{t=0}=[P,\tilde C]_{t=0},
\leqno{(3.3)}$$
where $\sigma_{P}$ is the Lichnerowicz coboundary (e.g., \cite{V3}).
Up to the sign, (3.3) is the coboundary defined in \cite{LMP} namely,
$$\sigma C=[\Lambda,C_{1}]-kE\wedge C_{1}-\Lambda\wedge C_{2}
\leqno{(3.4)}$$
$$-\frac{\partial}{\partial t}
\wedge([\Lambda,C_{2}]-(k-1)E\wedge C_{2}
+[E,C_{1}]).$$

We denote the Jacobi cohomology spaces by $H^{k}_{J}(M,\Lambda,E)$.
\proclaim 3.1 Remark. {\rm \cite{LMP}}. The anchor $\rho$ induces homomorphisms
$\rho^{\sharp}:H^{k}_{de\,R}(M)\rightarrow H^{k}_{J}(M,\Lambda,E)$ given
by
$$(\rho^{\sharp}\lambda)(s_{1},...,s_{k})=(-1)^{k}\lambda(\rho
s_{1},...,\rho s_{k})
\hspace{3mm}(\lambda\in\wedge^{k}M,\, s_{i}\in\Gamma J^{1}M).
\leqno{(3.5)}$$\par
Now, we need a recapitulation of several results of  \cite{Xu} and
\cite{ELW}.

For a Lie algebroid $A\rightarrow M$ with anchor $a$, an $A$-{\em connection}
$\nabla$ on a vector
bundle $E\rightarrow M$ consists of
{\em derivatives} $\nabla_{s}e\in\Gamma E$ ($s\in \Gamma
A,\;e\in\Gamma E$) which are ${\bf R}$-bilinear and satisfy
the conditions
$$\nabla_{fs}e=f\nabla_{s}e,\;\;
\nabla_{s}(fe)=(a(s)f)e+f\nabla_{s}e\hspace{3mm}(f\in C^{\infty}(M)).$$
For an $A$-connection, curvature may be defined as for usual connections.
Any flat $A$-connection $\nabla$ on $\wedge^{r}A$ ($r=rank\,A$)
produces a {\em Koszul operator} $D:\Gamma\wedge^{k}A
\rightarrow\Gamma\wedge^{k-1}A$,
locally given by
$$DU=(-1)^{r-k+1}[i(d\omega)\Omega+\sum_{h=1}^{r}\alpha^{h}
\wedge(i(\omega)\nabla_{s_{h}}\Omega],$$
where $\Omega\in\Gamma\wedge^{r}A$,
$\omega\in\Gamma\wedge^{r-k}A^{*}$ is such that $i(\omega)\Omega=U$,
$s_{h}$ is a local basis of $A$, and $\alpha^{h}$ is the dual cobasis of
$A^{*}$.
Moreover, $D$ is an exact generator of the
Gerstenhaber algebra of $A$, and every exact generator is defined by a
flat $A$-connection as above. The operator $D$ is a boundary, and yields
a corresponding homology, called the {\em homology of the Lie
algebroid $A$ with respect to the flat $A$-connection} $\nabla$,
$H_{k}(A,\nabla)$. For two flat connections $\nabla,\,\bar\nabla$
such that $D-\bar D=i(\alpha)$, where $\alpha=d_{A}f$ $(f\in C^{\infty}(M)
)$, one
has $H_{k}(M,\nabla)=H_{k}(M,\bar\nabla)$. If $\exists\Omega\in\Gamma
\wedge^{r}A$ which is nowhere zero, and
$\nabla\Omega=0$, one has the duality
$H_{k}(A,\nabla)=H^{r-k}(A)$, defined by sending $Q\in\Gamma\wedge^{k}A$
to $*_{\Omega}Q:=i(Q)\Omega$.

These results may be applied to the cotangent
Lie algebroid of an orientable Poisson manifold $(N^{n},Q)$. In this case,
the flat connection
$\nabla_{\theta}\Psi=\theta\wedge(di(Q)\Psi)$ ($\theta\in T^{*}N$,
$\Psi\in\wedge^{n}N$) 
precisely has the Koszul operator $\delta_{Q}$ and defines
the known Poisson homology $H_{k}(N,Q)$ (e.g., \cite{V3}). Finally (\cite{Xu},
Proposition 4.6 and Theorem 4.7), if $N$ has the volume form $\Omega$,
which defines a connection $\nabla_{0}$ by $\nabla_{0}\Omega=0$, and if
$W^{Q}$ is the {\em modular vector field} 
which acts on $f\in
C^{\infty}(M)$
according to the equation $L_{X^{Q}_{f}}\Omega=(W^{Q}f)\Omega$
($X^{Q}_{f}$ is
the Hamiltonian field of $f$) \cite{W},
one has $\delta_{Q}-D_{0}=i(W^{Q})$.
Accordingly, if the modular field $W^{Q}$ is Hamiltonian (i.e., $(N,Q)$
is a
{\em unimodular Poisson manifold}),
$H_{k}(N,Q)=H^{n-k}(T^{*}N)$.

The case of a general, possibly non orientable, Poisson manifold
is studied in \cite{ELW}. The expression of $\nabla_{\theta}\Psi$
above can be seen as the local equation of a connection on 
$\wedge^{n}T^{*}N$, and it still defines the Koszul operator
$\delta_{Q}$.
The general duality Theorem 4.5 of \cite{ELW} is
$$H_{k}(N,Q)=H^{n-k}(T^{*}N,\wedge^{n}T^{*}N),\leqno{(*)}$$
where the right hand side is the cohomology of the Lie algebroid $T^{*}N$
with values in the line bundle $\wedge^{n}T^{*}N$.
This means that the $k$-cocycles are spanned by cross sections
$V\otimes\Psi$, $V\in {\cal V}^{k}N$, $\Psi\in\Gamma\wedge^{n}T^{*}N$, and
the coboundary is given by
$$\partial(V\otimes\Psi)=[Q,V]\otimes\Psi+(-1)^{k}V\otimes\nabla\Psi.$$
Duality ($*$) is again defined 
by the isomorphism which sends $V\otimes\Psi$
to $i(V)\Psi$.

With this recapitulation finished, we apply the results to Jacobi
manifolds $(M^{m},\Lambda,E)$. Consider the Poisson manifold $(M\times{\bf
R},P)$ which we already used before.
Then $\delta_{P}$ is the Koszul operator of
the $(T^{*}M\times{\bf R})$-connection
$$\nabla_{\theta}\Psi=\theta\wedge(di(P)\Psi)\hspace{3mm}(\theta\in T^{*}
(M\times{\bf R}),\,\Psi\in\wedge^{m+1}(M\times{\bf R})).\leqno{(3.6)}$$
In particular, if we take
$$\theta=e^{t}(\alpha+fdt),\hspace{2mm}\Psi=e^{(m+1)t}\Phi\wedge dt
\hspace{3mm}(\alpha\in T^{*}M,\,\Phi\in\wedge^{m}M)\leqno{(3.7)}$$
(and use (2.9)) we get
$$\nabla_{\theta}\Psi=e^{(m+1)t}[fdi(E)\Phi-\alpha\wedge(di(\Lambda)\Phi
+mi(E)\Phi)]\wedge dt,\leqno{(3.8)}$$
and this formula may be seen as defining a $J^{1}M$-connection on
$\wedge^{m+1}J^{1}M$. Clearly, the Koszul operator of this connection must
be $\delta$ of (2.10). Therefore, we have
\proclaim 3.2 Proposition. The Jacobi homology of $(M,\Lambda,E)$ is
equal to the homology of the Lie algebroid $J^{1}M$ with respect to the
flat connection (3.8) i.e.,
$$H^{J}_{k}(M,\Lambda,E)=H_{k}(J^{1}M,\nabla).\leqno{(3.9)}$$ \par
Now, assume that $M$ has a volume form $\Phi\in\wedge^{m}M$. Then
$\Omega:=e^{(m+1)t}\Phi\wedge dt$ is a volume form on $M\times{\bf R}$, and
one has a connection $\nabla_{0}$ defined by
$\nabla_{0}\Omega=0$ with a Koszul operator $D_{0}$
such that $$\delta_{P}-D_{0}=i(W^{P}),\leqno{(3.10)}$$
where $W^{P}$ is the corresponding modular vector field i.e.,
$$L_{X^{P}_{\varphi}}\Omega=(W^{P}\varphi)\Omega
\hspace{3mm}(\varphi\in C^{\infty}(M\times{\bf R})).\leqno{(3.11)}$$

We need the interpretation of (3.10) at $t=0$. To get it, we take
local coordinates $(x^{i})$ on $M$, and compute the local components of
$W^{P}$ by using (3.11) for $\varphi=x^{i}$ and $\varphi=t$.
Generally,
we have
$$X^{P}_{\varphi}=i(d\varphi)P=e^{-t}(\sharp_{\Lambda}d\varphi+
\frac{\partial\varphi}{\partial t}E-(E\varphi)\frac{\partial}{\partial t}).
\leqno{(3.12)}$$
On the other hand, on $M$, let us define a vector field $V$
and a function $div_{\Phi}E$ by
$$L_{\sharp_{\Lambda}df}\Phi=(Vf)\Phi,\;L_{E}\Phi=(div_{\Phi}E)\Phi
\hspace{3mm}(f\in C^{\infty}(M). \leqno{(3.13)}$$
(The fact that $V$ is a derivation of $C^{\infty}(M)$ follows easily
from the skew symmetry of $\Lambda$.)
Then, the calculation of the local components of $W^{P}$ yield
$$W^{P}=e^{-t}[V-mE+(div_{\Phi}E)\frac{\partial}{\partial
t}].\leqno{(3.14)}$$
At $t=0$, (3.14) defines a section of $TM\oplus{\bf R}$ which we
denote by $V^{(\Lambda,E)}$, and call the {\em modular field} (not a
vector
field, of course) of the Jacobi manifold.

As in the Poisson case, if $\Phi\mapsto a\Phi$ ($a>0$),
$V^{(\Lambda,E)}\mapsto V^{(\Lambda,E)}+\sigma(\ln a)$ hence, what is
well
defined is the Jacobi cohomology class $[V^{(\Lambda,E)}]$, to be called
the
{\em modular class}.
If the modular class is zero $(M,\Lambda,E)$ is
a {\em unimodular Jacobi manifold}.

It is also possible to get the modular class $[V^{(\Lambda,E)}]$
from the general definition of the modular class of a Lie algebroid
\cite{ELW}. In the case of the algebroid $J^{1}M$, the definition of
\cite{ELW} means computing the expression
$${\cal E}:=(L^{J^{1}M}_{e^{t}(df+fdt)}[(e^{mt}\Phi)
\wedge(e^{t}dt)])\otimes\Phi$$
$$+(e^{(m+1)t}\Phi\wedge dt)\otimes(L_{\rho(e^{t}(df+fdt))}\Phi,$$
where $\rho$ is given by (2.6), and
$$L^{J^{1}M}_{e^{t}(df+fdt)}[(e^{mt}\Phi)
\wedge(e^{t}dt)]=\{e^{t}(df+fdt),(e^{mt}\Phi)
\wedge(e^{t}dt)\}_{P}.$$
If we decompose $(e^{mt}\Phi)=\wedge^{n}_{i=1}(e^{t}\varphi_{i})$,
$\varphi_{i}\in\wedge^{1}M$, the result of the required
computation turns out to be
$${\cal E}=(2(Vf)+2f(div_{\Phi}E)-Ef)(e^{(m+1)t}\Phi\wedge dt)
\otimes\Phi.$$
By comparing with (3.14), we see that the modular class in the sense of
\cite{ELW} is the Jacobi cohomology class of the cross section of
$TM\oplus {\bf R}$ defined by
$$A^{(\Lambda,E)}=2V^{(\Lambda,E)}-(2m+1)E.$$

With all this notation in place, the recalled results
of \cite{Xu}, Proposition 4.6, and \cite{ELW}, Theorem 4.5, yield
\proclaim 3.3 Proposition. If $(M,\Lambda,E)$ is a unimodular Jacobi
manifold one has duality between Jacobi homology and cohomology:
$$H^{J}_{k}(M,\Lambda,E)=H^{m-k+1}_{J}(M,\Lambda,E).\leqno{(3.15)}$$
If $(M,\Lambda,E)$ is an arbitrary Jacobi
manifold, one has the duality
$$H^{J}_{k}(M,\Lambda,E)=H^{m-k+1}_{J}(J^{1}M,
\wedge^{m+1}J^{1}M).\leqno{(3.15')}$$
\par
\noindent{\bf Proof.} 
The right hand side of (3.15$'$) is {\em Jacobi cohomology 
with values in} $\wedge^{m+1}J^{1}M$, similar to that in ($*$).
The homologies and cohomologies of (3.15) and (3.15$'$) are to be seen
as given by subcomplexes of $\oplus_{k}\wedge^{k}(M\times{\bf R})$,
$\oplus_{k}{\cal V}^{k}(M\times{\bf R})$
defined by (2.7) and (3.1). Then the result follows by the proofs of
the theorems of \cite{Xu}, \cite{ELW} quoted earlier, if we notice that
$$i[e^{-kt}(C_{1}+\frac{\partial}{\partial t}\wedge C_{2})]
(e^{(m+1)t}\Phi\wedge dt)$$
$$=e^{(m-k+1)t}[(-1)^{m}i(C_{2})\Phi+(i(C_{1})\Phi)\wedge dt].$$
The notation is that of (3.1) and (3.7). Q.e.d.

To get examples, let us look at the {\em transitive Jacobi manifolds}
\cite{GL}.

a). Let $M^{2n}$ be a locally or globally conformal symplectic manifold with the
$2$-form $\Omega=e^{\sigma_{\alpha}}\Omega_{\alpha}$, where
$\Omega_{\alpha}$ are
symplectic forms on the sets $U_{\alpha}$ of an open covering of $M$,
and $\sigma_{\alpha}\in
C^{\infty}(U_{\alpha})$. Then (e.g., \cite{V2}) $\{d\sigma_{\alpha}\}$ glue
up to a global closed $1$-form $\omega$, which is exact iff
$\exists\alpha$, $U_{\alpha}=M$, and $\sharp_{\Lambda}:=\flat_{\Omega}^{-1}$,
$E:=\sharp_{\Lambda}\omega$ define a Jacobi structure on $M$ \cite{GL}. It
follows easily that $L_{E}\Omega=0$ hence, $div_{\Omega^{n}}E=0$.
Furthermore,
$$L_{\sharp_{\Lambda}df}\Omega^{n}=-n(n-1)df\wedge\omega\wedge\Omega^{n-1}
.$$
Using the {\em Lepage decomposition theorem} (\cite{LM}, pg.46) we see that
$df\wedge\omega=\xi+\varphi\Omega$, where
$$\xi\wedge\Omega^{n-1}=0,\;\;\varphi=-\frac{1}{n}i(\Lambda)(df\wedge\omega)
=Ef.$$ 
Hence, $V=-n(n-1)E$, and
$V^{(\Lambda,E)}=-n(2n-1)E$. Then, for $f\in C^{\infty}(M)$,
(3.4) yields $\sigma f=\sharp_{\Lambda}df-(Ef)(\partial/\partial t)$, and
$\sigma f=E$ holds iff $\omega=df$.
Thus, (3.15) holds on  globally conformal symplectic manifolds. But,
(3.15) may not hold in the true locally conformal symplectic case. For
instance, it
follows from Corollary 3.15 of \cite{LMP0} that the result does not hold
on a Hopf manifold with its natural locally conformal K\"ahler structure.
(Private correspondence from J. C. Marrero.)

b). Let $M^{2n+1}$ be a contact manifold with the contact $1$-form $\theta$
such that $\Phi:=\theta\wedge(d\theta)^{n}$ is nowhere zero. Then $M$ has
the Reeb vector field $E$ where
$$i(E)\theta=1,\;\;i(E)d\theta=0,$$
and $\forall f\in C^{\infty}(M)$ there is a {\em Hamiltonian vector field}
$X^{\theta}_{f}$ such that
$$i(X^{\theta}_{f})\theta=f,\;\;i(X^{\theta}_{f})d\theta=-df+(Ef)\theta.$$
Furthermore, if $$\Lambda(df,dg):=d\theta(X^{\theta}_{f},X^{\theta}_{g})
\hspace{3mm}(f,g\in C^{\infty}(M)),$$
$(\Lambda,E)$ is a Jacobi structure \cite{GL}.

Now, let $(q^{i},p_{i},z)$ $(i=1,...,n)$ be local canonical coordinates,
such that $\theta=dz-\sum_{i}p_{i}dq^{i}$. Then
$$E=\frac{\partial}{\partial z},\;\Lambda=\sum_{i}\frac{\partial}{\partial
q^{i}}\wedge\frac{\partial}{\partial p_{i}}+\frac{\partial}{\partial z}
\wedge(\sum_{i}p_{i}\frac{\partial}{\partial p{_i}}).$$
This leads to $div_{\Phi}E=0$, $V^{(\Lambda,E)}=nE$, and it follows that
there is no $f\in C^{\infty}(M)$ satisfying $\sigma f=(nE\oplus 0)$.

We close this section by the remark that the identification of a manifold
$M$ with $M\times\{0\}\subseteq M\times {\bf R}$ leads to other interesting
structures too. For instance, if we define the spaces
$$\wedge_{0}^{k}M:=\{e^{t}(\xi_{1}+\xi_{2}\wedge dt)
\,/\,\xi_{1}\in\wedge^{k}M,\xi_{2}\in\wedge^{k-1}M\},$$
the triple $(\oplus_{k}
\wedge_{0}^{k}M,d,i(X+f(\partial/\partial t))$ is a {\em
Gelfand-Dorfman complex}
\cite{Dor}, and a Jacobi structure on $M$ is equivalent with
a {\em Hamiltonian structure} \cite{Dor} on this complex.

On the other hand, if we have a Jacobi manifold $(M,\Lambda,E)$, and put
$${\cal V}_{0}^{k}M:=\{e^{-(k-1)t}(Q_{1}+\frac{\partial}{\partial t}
\wedge Q_{2})\,/\,Q_{1}\in{\cal V}^{k}M,\,Q_{2}\in{\cal V}^{k-1}M\},$$
then $(\oplus_{k}{\cal V}_{0}^{k}M,[\:,\:],\sigma_{P})$ ($P$ is defined by
(2.2)) is a differential graded Lie algebra, the cohomology of which is
exactly the $1$-{\em differentiable Chevalley-Eilenberg cohomology}
$H^{k}_{1-dif}(M,\Lambda,E)$ of Lichnerowicz \cite{Lz}. In particular,
$H^{1}_{1-dif}(M,\Lambda,E)$ is the quotient of the space of conformal
Jacobi infinitesimal automorphisms by the space of the Jacobi Hamiltonian
vector fields \cite{Lz}.
\section{Lie bialgebroid structures on TM$\oplus${\bf R}}
In the Poisson case, $T^{*}M$ is a Lie bialgebroid
\cite{YK}, \cite{MX} with dual $TM$. This is
not true for $J^{1}M$ on Jacobi manifolds in spite of the fact
that $(J^{1}M)^{*}=TM\oplus{\bf R}$ has a natural Lie algebroid structure,
which extends the one of $TM$. Namely, if we
see ${\cal X}\in\Gamma(TM\oplus{\bf R})$ as a vector field of $M\times {\bf
R}$ given by
$${\cal X}=(X+f\frac{\partial}{\partial t})_{t=0}\hspace{3mm}
(X\in\Gamma TM,\,f\in C^{\infty}(M)),\leqno{(4.1)}$$
we have the Lie bracket
$$[{\cal X},{\cal Y}]_{0}:=
[X+f\frac{\partial}{\partial t},Y+g\frac{\partial}{\partial t}]=
[X,Y]+(Xg-Yf)\frac{\partial}{\partial t}, \leqno{(4.2)}$$
and the anchor map $a({\cal X}):=X$.
If we would be in the case of a Lie bialgebroid, the bracket
$\{f,g\}_{s}$
$:=<df,d_{*}g>$ $(f,g\in C^{\infty}(M))$,
where $d,\,d_{*}$ are the differentials of the Lie
algebroids $TM\oplus{\bf R}$ and $J^{1}M$, respectively, would be Poisson
\cite{YK}, \cite{MX}. This is not true since one gets
$\{f,g\}_{s}=\Lambda(df,dg)$.
\proclaim 4.1 Remark. The differential $\bar d$ defined by (2.11) is the
same as the differential $d$ of the Lie algebroid
$TM\oplus{\bf R}$ with the bracket (4.2).\par

In Poisson geometry, the cotangent Lie bialgebroid
structure is produced by a Poisson
bivector $\Pi$ of $TM$ i.e., $[\Pi,\Pi]=0$. It is natural to ask what is
the structure produced by a Poisson bivector $\Pi$ of $TM\oplus{\bf R}$. We
generalize this question a bit namely,
we fix a closed $2$-form $\Omega$ on $M$, and take the Lie bracket
$$[{\cal X},{\cal Y}]_{\Omega}:=[{\cal X},{\cal Y}]_{0}+\Omega(X,Y)
\frac{\partial}{\partial t}. \leqno{(4.2')}$$
The notation, and the anchor map $a$ are the same as for (4.2). It is known
that (4.2$'$) defines all the transitive Lie algebroid structures over $M$
such that the kernel of the anchor is a trivial line bundle, up to an
isomorphism \cite{Mk}. A Poisson
bivector $\Pi$ on $TM\oplus{\bf R}$ with the bracket (4.2$'$) will be called
an $\Omega$-{\em Poisson structure} on $M$.
\proclaim 4.2 Proposition. An $\Omega$-Poisson structure $\Pi$ on $M$ is
equivalent with a pair $(Q,E)$, where $Q$ is a Poisson bivector on $M$
(i.e., $[Q,Q]=0$), and $E$ is a vector field such that
$$L_{E}Q=\sharp_{Q}\Omega. \leqno{(4.3)}$$
\par \noindent{\bf Proof.}  Using the identification (4.1) of the 
cross sections of $TM\oplus{\bf R}$ with vector fields on $M\times{\bf R}$
for $t=0$, and
local coordinates $(x^{i})$ on $M$, we may
write
$$\Pi=Q+\frac{\partial}{\partial t}\wedge E=
\frac{1}{2}Q^{ij}(x)\frac{\partial}{\partial x^{i}}\wedge
\frac{\partial}{\partial x^{j}}+\frac{\partial}{\partial t}
\wedge(E^{k}(x)\frac{\partial}{\partial x^{k}}),\leqno{(4.4)}$$
where $Q$ is a bivector field on $M$, $E$ is a vector field, and the
Einstein
summation convention is used.

Now, $[\Pi,\Pi]_{\Omega}=0$ can be expressed by the known formula of the
Schouten-Nijenhuis bracket of decomposable multivectors (e.g., \cite{V3},
formula (1.12)), and (4.2$'$). The result is equivalent to $[Q,Q]=0$ and
(4.3). Q.e.d.
\proclaim 4.3 Corollary. If $(M,Q)$ is a Poisson manifold, $Q$ extends to a
$\Omega$-Poisson structure for every closed $2$-form $\Omega$,
where the de Rham class $[\Omega]$ has zero
$\sharp_{Q}$-image in the Poisson cohomology of $(M,Q)$, by taking $E$ such
that (4.3) holds. \par
This is just a reformulation of Proposition 4.2.

It is well known that a Poisson bivector on a Lie algebroid $A$ induces a
bracket on $\Gamma A^{*}$ such that $(A,A^{*})$ is a
triangular Lie bialgebroid
\cite{YK}, \cite{MX}. Namely, the Poisson bivector
$\Pi$  of (4.4) yields the following bracket
$$\{\alpha\oplus f,\beta\oplus g\}_{\Omega}:=
L^{\Omega}_{\sharp_{\Pi}(\alpha\oplus f)}(\beta\oplus g)
-L^{\Omega}_{\sharp_{\Pi}(\beta\oplus g)}(\alpha\oplus f)
-d_{\Omega}(\Pi(\alpha\oplus f,\beta\oplus g)),\leqno{(4.5)}$$
where $\alpha\oplus f,\beta\oplus g\in\Gamma J^{1}M$,
and the index $\Omega$ denotes the fact that the operators involved are
those of the Lie algebroid calculus of (4.2$'$).

To make this formula explicit, notice that
$$\sharp_{\Pi}(\alpha+fdt)=\sharp_{Q}\alpha+fE-\alpha(E)
\frac{\partial}{\partial t},\leqno{(4.6)}$$
whence $$\Pi(\alpha+fdt,\beta+gdt)=
Q(\alpha,\beta)+f\beta(E)-g\alpha(E).\leqno{(4.7)}$$
Then, by evaluation on a field of the form (4.1), and with (4.2$'$), we
obtain $$L_{\sharp_{\Pi}(\alpha+fdt)}^{\Omega}(\beta+gdt)=
L_{\sharp_{\Pi}(\alpha+fdt)}(\beta+gdt)-g(\flat_{\Omega}\sharp_{Q}
\alpha)-fgi(E)\Omega,\leqno{(4.8)}$$
where $\flat_{\Omega}X:=i(X)\Omega$.
As a consequence, (4.5) becomes
$$\{\alpha\oplus f,\beta\oplus g\}_{\Omega}:=
[\{\alpha,\beta\}_{Q}+f(L_{E}\beta+\flat_{\Omega}\sharp_{Q}\beta)
\leqno{(4.9)}$$
$$-g(L_{E}\alpha+\flat_{\Omega}\sharp_{Q}\alpha)]\oplus
[(\sharp_{Q}\alpha)g-(\sharp_{Q}\beta)f+f(Eg)-g(Ef)],$$
where $(Q,E)$ are associated with $\Pi$ as in Proposition 4.2.

The anchor map of the Lie algebroid $J^{1}M$ with (4.9)
is $\rho:=pr_{TM}\circ\sharp_{\Pi}$, and it is provided by (4.6).

In particular, Proposition 4.2 tells us that a pair $(Q,E)$ which consists
of a Poisson bivector $Q$ and an infinitesimal automorphism $E$ of $Q$, to
which we will refer as an {\em enriched Poisson structure}, provides a
Poisson bivector $\Pi$ on $TM\oplus{\bf R}$ with the bracket (4.2), and a
Lie bialgebroid
$(TM\oplus{\bf R},J^{1}M=T^{*}M\oplus{\bf R})$.

An example (suggested by \cite{Lh}) can be obtained as follows. Let
$(M,\Lambda,E)$ be a Jacobi manifold. A {\em time function} is a function
$\tau\in C^{\infty}(M)$ which has no critical points and satisfies $E\tau=1$.
If
such a function exists, $(\Lambda_{0}:=\Lambda-(\sharp_{\Lambda}dt)\wedge E,
E)$ is an enriched Poisson structure. Jacobi manifolds with time may be
seen as generalized phase spaces of time-dependent Hamiltonian systems.
Namely, if $H\in C^{\infty}(M)$ is the Hamiltonian function, the
trajectories of the system are the integral lines of the vector field
$X^{0}_{H}:=\sharp_{\Lambda_{0}}df+E$.

Let us briefly indicate the important objects associated with the 
Lie algebroids $TM\oplus{\bf R}$, defined by the bracket (4.2$'$),
and $J^{1}M$ with the bracket (4.9).

The cohomology of $TM\oplus{\bf R}$ is that of the cochain spaces
$$\wedge^{k}_{\Omega}M:=\{\lambda=\lambda_{1}+\lambda_{2}\wedge dt
\;/\;\lambda_{1}\in \wedge^{k}M,\lambda_{2}\in\wedge^{k-1}M\}
\leqno{(4.10)}$$
with the corresponding coboundary, say $d_{\Omega}$. A straightforward
evaluation of $d_{\Omega}\lambda$ on arguments
$X_{i}+f_{i}(\partial/\partial t)$, in accordance with the Lie algebroid
calculus \cite{Mk}, yields the formula
$$d_{\Omega}\lambda=d\lambda-(-1)^{k}\Omega\wedge\lambda_{2}.
\leqno{(4.11)}$$

The {\em Poisson cohomology} of $TM\oplus{\bf R}$ above i.e., the
cohomology of the Lie algebroid $J^{1}M$ with (4.9), can be seen (with
(4.1)) as having the cocycle spaces
$${\cal C}^{k}(M):=\{C=C_{1}+\frac{\partial}{\partial t}\wedge C_{2}\;
/\;C_{1}\in{\cal V}^{k}M,C_{2}\in {\cal V}^{k-1}M\},\leqno{(4.12)}$$
and the coboundary $\partial C=[\Pi,C]_{\Omega}$, with $\Pi$ of (4.4) and
the $\Omega$-Schouten-Nijenhuis bracket. In order to write down a concrete
expression of this coboundary, we define an operation
$U\wedge_{\Omega}V\in{\cal V}^{k+h-2}$, for $U\in {\cal V}^{k}M$,
$V\in{\cal V}^{h}M$, by the formula
$$U\wedge_{\Omega}V(\alpha_{1},...,\alpha_{k+h-2})=\frac{1}
{(k-1)!(h-1)!}\sum_{\sigma\in
S_{k+h-2}}[(sign\,\sigma)\leqno{(4.13)}$$
$$\cdot\sum_{i=1}^{m}U(\epsilon^{i},
\alpha_{\sigma_{1}},...,\alpha_{\sigma_{k-1}})
V(\flat_{\Omega}e_{i},\alpha_{\sigma_{k}},...,\alpha_{\sigma_{k+h-2}})],
$$
where $S$ is the symmetric group, $e_{i}$ is a local tangent basis
of $M$, and $\epsilon^{i}$ is the corresponding dual cobasis.
If $U,V$ are vector fields, $U\wedge_{\Omega}V=\Omega(U,V)$.
By computing for decomposable multivectors $C_{1},C_{2}$, we get
$$\partial C=[Q,C_{1}]+\frac{\partial}{\partial t}\wedge
([Q,C_{2}]+Q\wedge_{\Omega} C_{1}-L_{E}C_{1}),\leqno{(4.14)}$$
where the brackets are the usual Schouten-Nijenhuis brackets on $M$.

Furthermore, the exact generator of the BV-algebra of the Lie algebroid
$J^{1}M$ is $\delta_{\Omega}:=[i(\Pi),d_{\Omega}]$, and using (4.11) we get
$$\delta_{\Omega}(\lambda_{1}+\lambda_{2}\wedge dt)=
\delta_{Q}\lambda_{1}+(-1)^{k-1}([i(Q),e(\Omega)]\lambda_{2}\leqno{(4.15)}$$
$$-di(E)\lambda_{2})
+(\delta_{Q}\lambda_{2})\wedge dt,$$
where $e(\Omega)$ is exterior product and $[\;,\;]$ is the commutant of
the operators.

Finally, let us discuss the {\em modular class} of the Lie algebroid
($J^{1}M$, (4.9)). For simplicity, we assume the manifold $M$ orientable,
with the volume form $\Phi\in\Gamma\wedge^{m}M$. In the non orientable
case, the same computations hold if $\Phi$ is replaced by a density 
$\Phi\in\Gamma|\wedge^{m}M|$. 
Again, we denote by $W^{Q}$ the modular vector
field defined by $L_{\sharp_{Q}df}\Phi=(W^{Q}f)\Phi$ (see Section 3).

There are two natural possibilities to define a modular class for 
the algebroid $J^{1}M$. One is by computing the Lie derivative:
$$L_{\sharp_{\Pi}(\alpha+fdt)}(\Phi\wedge dt)=
[L_{\sharp_{Q}\alpha}\Phi+fL_{E}\Phi+df\wedge i(E)\Phi]
\wedge dt.\leqno{(4.16)}$$
This result is obtained if the computation is done after $\Phi$
is decomposed into a product of $m$ $1$-forms, and by using (4.6). Since
$i(E)(df\wedge\Phi)=0$, the last term in (4.16) is $(Ef)\Phi$, and if we
also use (3.13), (4.16) yields
$$L_{\sharp_{\Pi}(df+fdt)}(\Phi\wedge dt)=
(W^{Q}f+fdiv_{\Phi}E+Ef)(\Phi\wedge dt).\leqno{(4.17)}$$
Therefore, we get the {\em modular field}
$$W^{\Pi}:=W^{Q}+E+(div_{\Phi}E)\frac{\partial}{\partial
t}.\leqno{(4.18)}$$
If $\Phi$ is replaced by $h\Phi$ $(h\in C^{\infty}(M))$, it follows easily
that the $\Pi$-Poisson cohomology class $[W^{\Pi}]$ is preserved. This
will be the {\em modular class}.

The second possibility is to apply the general definition of \cite{ELW}.
Similar to what we had for Jacobi manifolds in Section 3,
this asks us to compute the flat connection $D$ on $(\wedge^{m+1}J^{1}M)
\otimes(\wedge^{m}T^{*}M)$ given by
$$D_{(df+fdt)}[(\Phi\wedge dt)\otimes\Phi]=
L^{J^{1}M}_{(df+fdt)}(\Phi\wedge dt)\otimes\Phi+
(\Phi\wedge dt)\otimes (L_{\rho(df+fdt)}\Phi).\leqno{(4.19)}$$

From Lie algebroid calculus, we know that
$$L^{J^{1}M}_{(df+fdt)}(\Phi\wedge dt)=\{df+fdt,\Phi\wedge dt\}_{\Omega},
\leqno{(4.20)}$$
where the bracket is the Schouten-Nijenhuis extension of (4.9).
If we look at a decomposition $\Phi=\varphi_{1}\wedge...\wedge
\varphi_{n}$ $(\varphi_{i}\in\wedge^{1}M)$, (4.9) yields
$$\{df+fdt,dt\}_{\Omega}=0,$$
$$\{df+fdt,\varphi_{i}\}_{\Omega}=L_{\rho(df+fdt)}\varphi_{i}
-\varphi_{i}(E)df+f\flat_{\Omega}\sharp_{Q}\beta-(\sharp_{Q}\beta(f))dt,$$
and we get
$$L^{J^{1}M}_{(df+fdt)}(\Phi\wedge dt)=[L_{(df+fdt)}\Phi-(Ef)\Phi
+ftr(\flat_{\Omega}\circ\sharp_{Q})\Phi]\wedge dt.\leqno{(4.21)}$$

But, we also have
$$L_{\rho(df+fdt)}\Phi=L_{\sharp_{Q}df+fE}\Phi=[W^{Q}f+fdiv_{\Phi}E
+Ef]\Phi.\leqno{(4.22)}$$
By inserting (4.21), (4.22) into (4.19), we get another {\em modular
field} namely,
$$A_{\Omega}:=(2W^{Q}+E)\oplus(div_{\Phi}E+tr(\flat_{\Omega}\circ\sharp_{Q}))
=(2W^{\Pi}-E)\oplus tr(\flat_{\Omega}\circ\sharp_{Q}).\leqno{(4.23)}$$
From the general results of \cite{ELW}, it is known that the $\Pi$-Poisson
cohomology class of this field is independent of the choice of $\Phi$, and
it is a {\em modular class} of $J^{1}M$.

As for the modular class of $TM\oplus{\bf R}$ with the bracket (4.2$'$),
it vanishes for reasons similar to those for the class of the tangent
groupoid $TM$ \cite{ELW}.

We finish by another interpretation of the enriched Poisson
structures.
If ${\cal F}$, be an arbitrary associative, commutative, real algebra, we
may say that $f:M\rightarrow{\cal F}$ is differentiable if for any ${\bf
R}$-linear mapping $\phi:{\cal F}\rightarrow{\bf R}$, $\phi\circ f\in
C^{\infty}(M)$. Furthermore, an ${\bf R}$-linear operator $v_{x}$ which
acts on germs of ${\cal F}$-valued differentiable functions at $x\in M$,
and satisfies the Leibniz rule will be an ${\cal F}$-{\em tangent vector}
of $M$ at $x$. Then, we have natural definitions of tangent spaces
$T_{x}(M,{\cal F})$,
differentiable ${\cal F}$-vector fields, etc. \cite{V1}.
A bracket $\{\;,\;\}$ which makes $C^{\infty}(M,{\cal F})$ into a Poisson
algebra will be called an ${\cal F}$-{\em Poisson structure} on $M$.

Now, take ${\cal F}$ to be the {\em Studi algebra of parabolic dual numbers}
${\cal S}:={\bf R}[s\,/\,s^{2}=0]$, where $s$ is the generator. An ${\cal
S}$-Poisson structure $\Pi$ in the above mentioned sense will be called a
{\em Studi-Poisson structure}. The restriction of $\Pi$ to real valued
functions is a Poisson bivector $Q$ on $M$, and the Jacobi identity shows
that the Hamiltonian vector field $X^{\Pi}_{s}$
of the constant function $s$ is an infinitesimal
automorphism $E$ of $Q$. Conversely, the pair $(Q,E)$ defines the
Studi-Poisson bracket
$$\{f_{0}+f_{1}s,g_{0}+g_{1}s\}:=\{f_{0},g_{0}\}_{Q}+f_{1}(Eg_{0})-
g_{1}(Ef_{0})$$
$$+s(\{f_{0},g_{1}\}_{Q}+\{f_{1},g_{0}\}_{Q}+f_{1}(Eg_{1})-g_{1}(Ef_{1})).$$

Notice that we cannot say that $v_{x}s=0$ for all $v_{x}\in T_{x}(M,{\cal
F})$ since $v_{x}$ was linear over ${\bf R}$ only.
 \vspace*{1cm}
{\small Department of Mathematics, \\}
{\small University of Haifa, Israel.\\}
{\small E-mail: vaisman@math.haifa.ac.il}
\end{document}